\documentclass[conference]{IEEEtran}
\IEEEoverridecommandlockouts
\usepackage{cite}
\usepackage{amsmath,amssymb,amsfonts}
\usepackage{algorithmic}
\usepackage{graphicx}
\usepackage{textcomp}
\usepackage{xcolor}

\usepackage{graphicx}
\graphicspath{{./figs/}}
\usepackage{epstopdf}
\usepackage{url}
\usepackage{hyperref}

\usepackage[utf8]{inputenc}

\newcommand{\vect}[1]{\boldsymbol{\mathbf{#1}}}
\usepackage{mathtools}

\def\BibTeX{{\rm B\kern-.05em{\sc i\kern-.025em b}\kern-.08em
    T\kern-.1667em\lower.7ex\hbox{E}\kern-.125emX}}
\begin{document}

\title{Data-driven identification of vehicle dynamics using Koopman operator
}


\author{\IEEEauthorblockN{Vít Cibulka\IEEEauthorrefmark{1},
Tomáš Haniš\IEEEauthorrefmark{2} and
Martin Hromčík\IEEEauthorrefmark{3}}
\IEEEauthorblockA{Department of Control Engineering,
Faculty of Electrical Engineering, Czech Technical University in Prague\\ Email:
\IEEEauthorrefmark{1}cibulka.vit@fel.cvut.cz
\IEEEauthorrefmark{2}hanis.tomas@fel.cvut.cz
\IEEEauthorrefmark{3}hromcik.martin@fel.cvut.cz}}


\maketitle

\begin{abstract}
This paper presents the results of identification of vehicle dynamics using the Koopman operator. The basic idea is to transform the state space of a nonlinear system (a car in our case) to a higher-dimensional space, using so-called basis functions, where the system dynamics is linear.
The selection of basis functions is crucial and there is no general approach on how to select them, this paper gives some discussion on this topic.
Two distinct approaches for selecting the basis functions are presented. The first approach, based on Extended Dynamic Mode Decomposition, relies heavily on expert basis selection and is completely data-driven.
The second approach utilizes the knowledge of the nonlinear dynamics, which is used to construct eigenfunctions of the Koopman operator which are known by definition to evolve linearly along the nonlinear system trajectory. The eigenfunctions are then used as basis functions for prediction.
Each approach is presented with a numerical example and discussion on the feasibility of the approach for a nonlinear vehicle system.
\end{abstract}


\section{Introduction}


 Car is a highly nonlinear dynamical system which has been used for decades without any sort of sophisticated automatic control system.
Our goal is to create full-time-full-authority control system, which means that the driver will rather set the reference for the car, and the ultimate control over the car inputs (steering angle, motor torques) will lie with the control system, protecting the driver from himself.
However, control of nonlinear systems is not an easy task. 

The Koopman operator is becoming an increasingly popular tool for nonlinear control. The Koopman operator is a linear operator which can in theory approximate any nonlinear system. The idea lies in \textit{lifting} the state space of the nonlinear system to a higher-dimensional space by means of a nonlinear transformation. The lifted system should then be evolving linearly along with the original nonlinear system. The lifted state space is then transformed back with a linear transformation to the original state space and the nonlinear states (or system outputs) are recovered.

Such an approach is suitable for any linear control-design method such as LQG, MPC and $H_{\infty}$.

In this paper, we present the results of application of the Koopman operator framework to vehicle system in an attempt to create a linear representation of a vehicle and use it to control the vehicle optimally with linear MPC.

\section{Singletrack model}
\begin{figure}[htbp]
\centerline{\includegraphics[width=0.5\textwidth]{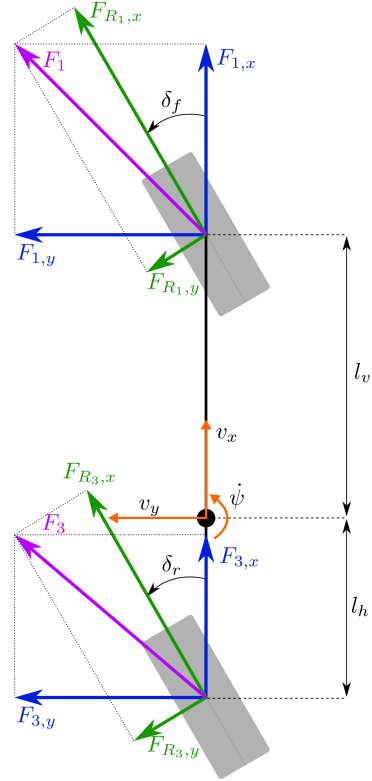}}
\caption{The singletrack model. Forces $F_{R_2}$ and $F_{R_4}$ are not depicted in the figure because in a general case with symmetric tires $F_{R_2} = F_{R_1}$ and $F_{R_4} = F_{R_3}$.}
\label{fig:singletrack}
\end{figure}
The singletrack model is depicted in Fig.\ref{fig:singletrack}.
The model was derived from 16-state twin-track model described in \cite{Schramm2014}.
The states and parameters of this model are a subset of states and parameters of the full twin-track model.

State vector of the singletrack model is 
$
\begin{bmatrix}
v_x (m/s)&
v_y (m/s)&
\dot{\psi} (rad/s)&
\dot{\rho}_f (rad/s)&
\dot{\rho}_r (rad/s)
\end{bmatrix}^T
$, where $v_x$ is longitudinal velocity, $v_y$ lateral velocity, $\dot{\psi}$   yawrate and $\dot{\rho}_{f/r}$ are front/rear wheel angular rates.

This model has 4 wheels, with two wheel always being in the same place. This allows for usage of asymmetric tire models, reduces the work needed to transition to twin-track model and is less error-prone than the standard approach (with only two tires, the user has to remember that each tire should generate twice as much force).

The vehicle body is modeled as a rigid body using Newton-Euler equations.
\begin{equation}
m_v ( 
\begin{bmatrix}
\dot{v}_x \\
\dot{v}_y \\
\end{bmatrix} + 
\dot{\psi}
\begin{bmatrix}
 -v_y \\
 v_x \\
 \end{bmatrix} )
 =
 \sum_{i=1}^{4} 
 \begin{bmatrix}
 \mathbf{F}_{i,x} \\
 \mathbf{F}_{i,y} \\
 \end{bmatrix}
 -\frac{1}{2}c_w \rho A_w \sqrt{v_x^2 + v_y^2}
 \begin{bmatrix}
 v_x\\v_y
 \end{bmatrix}
\end{equation}
\begin{equation}
J_{zz} \ddot{\psi}
=\sum_{i=1}^{4} \vect{r}_i \mathbf{F}_{i,y}
\label{eq:yawr}
\end{equation}
\begin{equation}
J_{R_i} \ddot{\rho}_{R_i} = M_{a,R_i} - M_{b,R_i} sign(\dot{\rho}_{R_i}) - rF_{R_i,x} ,\forall i = 1,3
\label{eq:torque}
\end{equation}
Where
\begin{equation}
\vect{r} = 
\begin{bmatrix}
 \vect{r_1}\\ \vect{r_2}\\\vect{r_3}\\\vect{r_4}
 \end{bmatrix}^T
  = 
  \begin{bmatrix}
  \begin{bmatrix}
  l_v\\0\\0
  \end{bmatrix},
  \begin{bmatrix}
  l_v\\0\\0
  \end{bmatrix},
  \begin{bmatrix}
  -l_h\\0\\0
  \end{bmatrix},
  \begin{bmatrix}
  -l_h\\0\\0
  \end{bmatrix}
  \end{bmatrix}
\end{equation} is the vector describing position of each wheel with respect to the center of gravity. The wheels are numbered in this order: front-left, front-right, rear-left, rear-right.
$ m_v $ is the vehicle mass,  
\(\mathbf{F}_{i,x/y}\) is a force acting on i-th wheel along x/y axis in body-fixed coordinates. $F_{R_i,x}$ is a force acting along x axis in wheel coordinate system (direct output of the tire model).
 The term \( -\frac{1}{2}c_w \rho A_w \sqrt{v_x^2 + v_y^2}
 \begin{bmatrix}
 v_x\\v_y
 \end{bmatrix}\)
 is an approximation of air-resistance, \(c_w\) is a drag coefficient, \(\rho\) is air density and \(A_w\) is the total surface exposed to the air flow.
\(J_{zz}\) is the vehicle inertia about z-axis.
\(J_{R_i}\) the wheel inertia about y-axis.
Inputs are wheel torques $M_{a,R_i}$ (throttle), $M_{b,R_i}$ (break) and steering angles $\delta_{f/r}$.

Since this is a 4-wheel singletrack model, the following holds:
\begin{equation}
\dot{\rho}_{R_1} = \dot{\rho}_{R_2}
\end{equation}
\begin{equation}
\dot{\rho}_{R_3} = \dot{\rho}_{R_4}.
\end{equation}

The forces \(  \begin{bmatrix}
 F_{R_i,x} \\
 F_{R_i,y} \\
 \end{bmatrix} \)
 are calculated using the ``Pacejka magic formula'' \cite{Pacejka2012}
 \begin{equation}
F = D \cos(C \arctan(Bx - E(Bx - \arctan(Bx)))).
\label{pacejka}\end{equation}
The same formula can be used for calculating $F_x$ (tire longitudinal force) and $F_y$ (tire lateral force) with a different set of parameters for each.
The argument $x$ can be either sideslip angle $\alpha$ or longitudinal slip $\kappa$ (see \cite{Pacejka2012}) for calculating $F_y$ or $F_x$ respectively.
The parameters B,C,D and E are usually time-dependent. This work uses the Pacejka tire model \cite{Pacejka2012} with coefficients from the \textit{Automotive challenge 2018} organized by Rimac Automobili.

The transformation of tire forces from wheel-coordinate system to car coordinate system is done as follows

\begin{equation}
\begin{bmatrix}
 F_{i,x} \\
 F_{i,y} \\
 \end{bmatrix}
=
\begin{bmatrix}
\cos(\delta_i) & -\sin(\delta_i) \\
\sin(\delta_i) & \cos(\delta_i)
\end{bmatrix}
\begin{bmatrix}
 F_{R_i,x} \\
 F_{R_i,y} \\
 \end{bmatrix}
\end{equation}

\subsection{3 state singletrack}
\label{sec:3state}
To further simplify the model, equation \eqref{eq:torque} can be omitted, resulting in a model with only 3 states: $
\begin{bmatrix}
v_x &
v_y &
\dot{\psi}
\end{bmatrix}^T
$.
The inputs are then longitudinal slips (which were previously derived from \eqref{eq:torque}) and steering angles.
This model will be used for the Eigenfunction approach in \ref{sec:eigenfun}.
\subsection{3 state singletrack without tire model}
\label{sec:simplified}
To simplify the model even more, one can omit the tire model and use the tire forces $\begin{bmatrix}
 F_{i,x} \\
 F_{i,y} \\
 \end{bmatrix}$ as input, assuming the existence of a higher level control system controlling the tire forces and thus securing the assumption that the car can be controlled directly by force reference.

This model will be used for the Extended Dynamic Mode Decomposition (EDMD) \cite{koop_edmd} approach in \ref{sec:edmd}.

\section{The Koopman Operator}
\subsection{Extended Dynamic Mode Decomposition approach}
\label{sec:edmd}
\subsubsection{Framework description}
The Koopman operator framework for controlled systems, as described in \cite{Korda_koop} will be reviewed now.
Let us assume uncontrolled discrete nonlinear dynamical system with state $x$ and dynamics $f(\cdot)$
\begin{align}
\begin{split}
x^+ &= f(x) \\
y &= g(x)
\label{eq:nldyn}
\end{split}
\end{align}
with $x$ being the current state and $x^+$ the next state.

The Koopman operator $\mathcal{K} : \mathcal{F} \rightarrow \mathcal{F}$ is defined as 
\begin{equation}
(\mathcal{K} \psi)(x) = \psi(f(x))
\end{equation}
for each basis function $\psi : \mathbb{R}^N \rightarrow \mathbb{R}$ where $\mathcal{F}$ is the space of basis functions.

$\mathcal{K}$ is infinite-dimensional linear operator, which can be approximated by EDMD. The approximation is done by solving optimization problem

\begin{equation}
\min_A \sum_{j=1}^{K} || \boldsymbol{\psi}(x_j^+) - A\boldsymbol{\psi}(x_j) ||_2^2
\label{eq:opt1}
\end{equation}
where $\boldsymbol{\psi} = \begin{bmatrix}
\psi_1(x) & \psi_2(x) & \ldots \psi_{N_\psi}(x)
\end{bmatrix}^T $. The states $x_j$ and $x_j^+$ are obtained by simulating the nonlinear system model, $K$ is the cardinality of the simulated dataset $\mathcal{D} = \begin{bmatrix}
x_j & x_j^+ & y_j\end{bmatrix}_{j=1}^K$. Note that the states $x_j$ do not have come from some trajectory of the system. When the system model is available, it is sufficient to sample the state space with $x_j$ and perform a one-step simulation to obtain the $x_j^+$. 

The matrix $A$ now defines a discrete linear system which approximates the nonlinear dynamics from \eqref{eq:nldyn}
\begin{align}
\begin{split}
z^+ &= Az \\
\hat{y} &= Cz
\end{split}
\end{align}
where $z = \boldsymbol{\psi}(x)$ and $\hat{y}$ is the output estimate.
The matrix $C$ is obtained in a similar fashion
\begin{equation}
\min_C \sum_{j=1}^{K} || y_j - C\boldsymbol{\psi}(x_j) ||_2^2
\label{eq:C}
\end{equation}
For controlled system
\begin{align}
\begin{split}
x^+ &= f(x,u) \\
y &= g(x)
\end{split}
\end{align}
the approach is similar. The optimization problem \eqref{eq:opt1} changes to 
\begin{equation}
\min_{A,B} \sum_{j=1}^{K} || \boldsymbol{\psi}(x_j^+) - (A\boldsymbol{\psi}(x_j) + Bu_j ) ||_2^2
\end{equation} while \eqref{eq:C} stays the same.
The Koopman operator $\mathcal{K}$ is approximated both by $A$ and $B$ and acts as a one-step predictor.
Note that the learning dataset now consists of $\mathcal{D} = \begin{bmatrix}
x_j & x_j^+ &y_j& u_j\end{bmatrix}_{j=1}^K$.
See \cite{Korda_koop} for more information.

This approach will be showed on the model which does not include the tire nonlinearities, described in \ref{sec:simplified}, because the selection of $\boldsymbol{\psi}$ and $\mathcal{D}$ is much more difficult when the tire nonlinearities are present.


\subsubsection{Data selection}
\label{sec:dataset1}
The learning dataset was selected as follows.
The velocities \(v_x\) and $v_y$ were uniformly sampled in the interval $ \begin{bmatrix}
-30 & 30
\end{bmatrix} m/s$,  yawrate was uniformly sampled in interval $ \begin{bmatrix}
-10 & 10
\end{bmatrix} rad/s$. The number of samples was 15 for each state, resulting in 3 \textit{1-by-15} vectors. 

The input forces to the model were sampled uniformly in the interval $ \begin{bmatrix}
-F_m & F_m
\end{bmatrix}N$ with $F_m \in \{1,5,10,100\}$, each interval with 15 samples, resulting in 4 \textit{4-by-15} vectors.
The columns of sampled vectors were combined into all possible combinations, resulting in $4\cdotp15^4$ values in the dataset $\mathcal{D}$.

Note that the tire force of a regular car is usually in order of thousands of Newtons. However, the results were much better with relatively small inputs forces. Large forces were too dominant and their influence on the system overshadowed the system dynamics and since there are no tire nonlinearities, the system can be sufficiently excited with small input forces.

\subsubsection{Basis selection}
\begin{figure}[htbp]
\centerline{\includegraphics[width=0.5\textwidth]{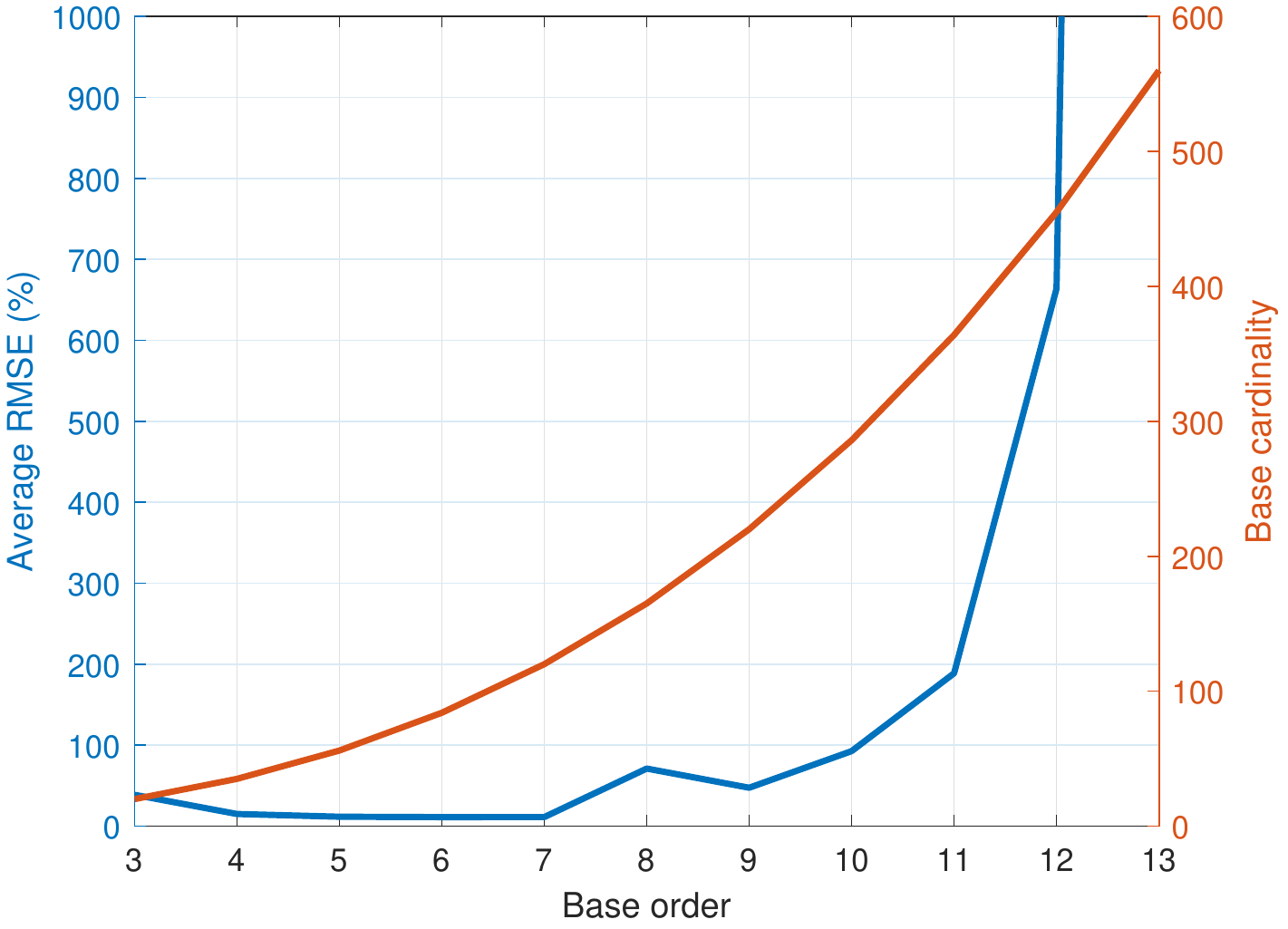}}
\caption{Comparison of polynomial bases with different orders. It is clear from the figure that large set of basis functions does not imply better prediction error. The best results were for polynomial bases of orders 6 and 7 with $RMSE = 11.1\%$.}
\label{fig:bases}
\end{figure}

There is unfortunately no general recommendation to determine which basis functions to choose for a given system. Some systems show satisfactory results with thin plate spline radial basis functions (see \cite{Korda_eigen}).
In our case, polynomial basis functions yielded much better results. A polynomial basis consists of monomials of the elements of the state vector. 
\begin{equation}
B_{k} = \{v_x^a \cdotp v_y^b \cdotp \dot{\psi}^c | a,b,c \in \{0,1,2,\ldots k\}\}
\end{equation}
where $k$ is the order of the basis $B_k$.
Bases with various orders were compared, the comparison was made using averaged root mean square error (RMSE)
\begin{equation}
RMSE = 100\frac{\sqrt{\sum_{k}{ || x_{koop}(kT_s) - x_{real}(kT_s)  ||_2^2}} }{\sqrt{ || \sum_{k} x_{real}(kT_s) ||_2^2}}
\label{eq:rmse}
\end{equation}
calculated over 3375 trajectories with the same distribution of initial conditions as the learning dataset $\mathcal{D}$. The length of each trajectory was 30 samples, with $T_s = 0.01$. The results can be seen in Fig.\ref{fig:bases}.

One might think the error would decrease with adding more basis functions but error grows exponentially as the order increases, which means that one cannot simply add more basis functions and expect the error to diminish.
When considering linear Koopman system
\begin{equation}
\boldsymbol{\dot{\psi}}(x) = A \boldsymbol{\psi}(x),
\end{equation}
the more basis functions are in $\boldsymbol{\psi}$, the more derivatives $\boldsymbol{\dot{\psi}}$ have to be expressed by $\boldsymbol{\psi}$. The derivatives of the polynomial basis grow in complexity with increasing order, which makes it harder to express them, so the result makes sense, although it might seem counter-intuitive at first glance.

\subsubsection*{Results}
\begin{figure*}[htbp]
\centerline{\includegraphics[width=\textwidth]{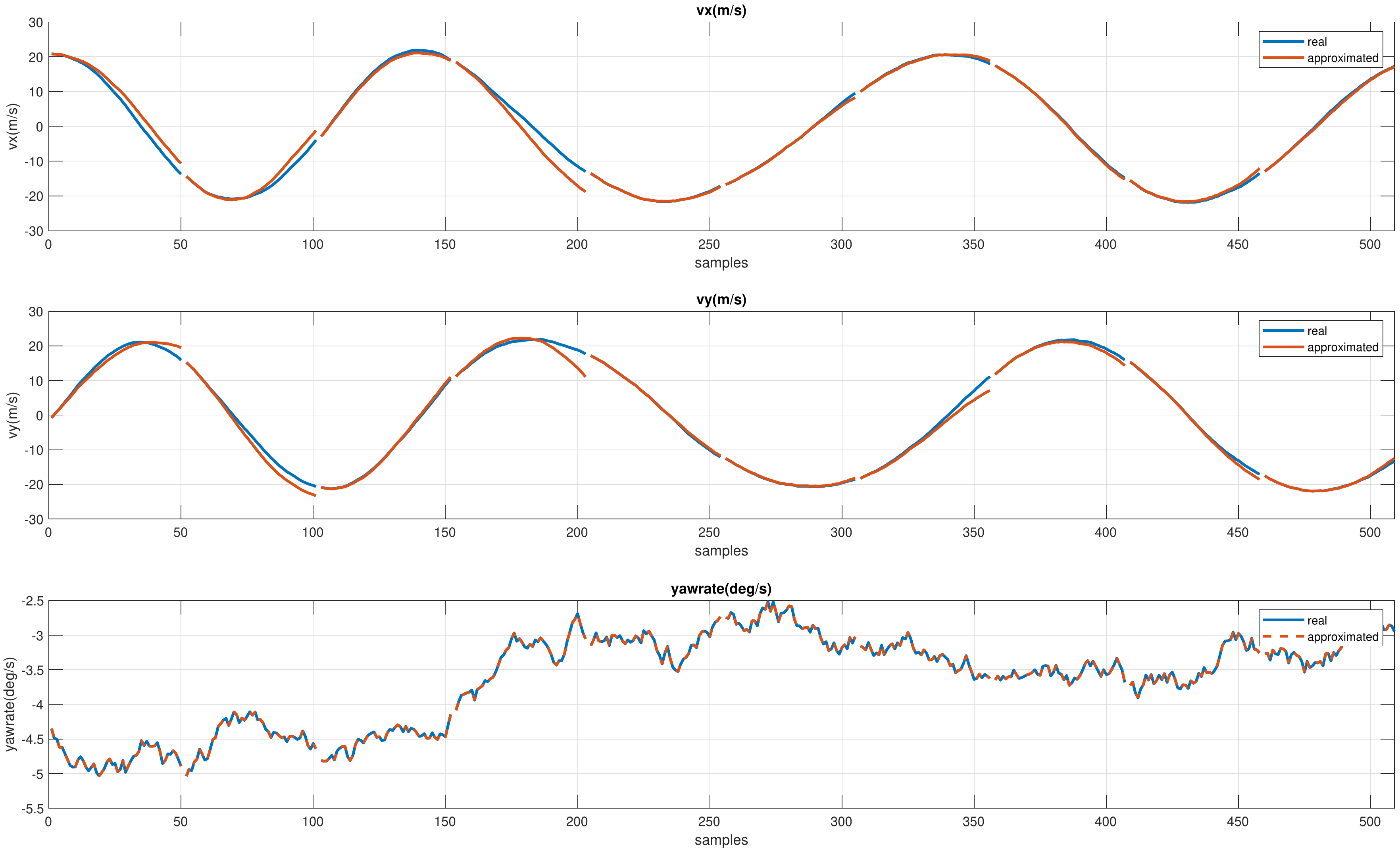}}
\caption{This figure shows the comparison between the real nonlinear system and linear system approximated by the Koopman operator with polynomial basis of order 7. The is input random with uniform distribution in the interval $[-10^4,10^4]$. The initial conditions were chosen randomly as $[v_x,v_y,\dot{\psi}] = 
  [20.8533
   ,-0.7222
   ,-4.3479]$.
The prediction is restarted after every 50 samples, the restart is visualized by segmenting the whole trajectory into pieces, each 50 samples long.}
\label{fig:koop1}
\end{figure*}

In our experiments, a polynomial basis with order 7 was used. The Koopman operator was approximated using the dataset described above.
The results can be seen in Fig. \ref{fig:koop1}. 
The yawrate is tracked without error because it is simply an integral of the input force as seen in \eqref{eq:yawr}.

\subsection{Basis functions from data}
\label{sec:eigenfun}
To deal with the problems of basis function selection, the approach from \cite{Korda_eigen} can be used. The approach consists of creating the eigenfunctions of the Koopman operator from data and using them as basis functions $\boldsymbol{\psi}$. 
The matrices $A$ and $B$ are calculated separately, contrary to the previous approach, which allows for optimizing the resulting system for a multiple-step prediction as opposed to one-step prediction of the previous approach. More in \cite{Korda_eigen}.

The basic idea of an eigenfunction will now be established.
Let us consider uncontrolled nonlinear continuous system
\begin{equation}
\frac{d}{dt}x = f(x)
\label{eq:koop2}
\end{equation}
Starting in $x_0$, the system will get to state $x_t$ in time $t$.
The state vector $x_t$ will be transformed with a function $\psi(x_t)$, defined as
\begin{equation}
\psi(x_t) = \psi(x_t)_{\lambda,g} = e^{\lambda t}g(x_0)
\label{eq:basis}
\end{equation}
for arbitrary eigenvalue $\lambda$ and function $g : \mathbb{R}^n \rightarrow \mathbb{R}$.
The time derivative of $\psi(x_t)$ is
\begin{equation}
\frac{d}{dt} \psi(x_t) = \frac{d}{dt} e^{\lambda t}g(x_0) = \lambda e^{\lambda t} g(x_0) = \lambda \psi(x_t) 
\end{equation}
We see that $\psi(x_t)$ evolves linearly with the system \eqref{eq:koop2}. For $\boldsymbol{\psi}(x_t) = \begin{bmatrix}
\psi_1(x_t) & \psi_2(x_t) & \ldots \psi_{N_\psi}(x_t)
\end{bmatrix}^T $ we would get
\begin{equation}
\frac{d}{dt}
\begin{bmatrix}
\psi_1(x_t) \\
 \psi_2(x_t) \\
  \vdots \\
  \psi_{N_\psi}(x_t)
\end{bmatrix}
 =
 \begin{bmatrix}
 \lambda_1 \\
 &\lambda_2\\
 &&\ddots\\
 &&&\lambda_{N_\psi}
 \end{bmatrix}
 \begin{bmatrix}
\psi_1(x_t) \\
 \psi_2(x_t) \\
  \vdots \\
  \psi_{N_\psi}(x_t)
\end{bmatrix}
\end{equation}
\begin{equation}
\frac{d}{dt} z = A z
\end{equation}
for $z = \boldsymbol{\psi}(x)$.
Choosing \eqref{eq:basis} as basis functions immediately yields the diagonal $A$ matrix. Notice that in this case, the Koopman system is derived as continuous system, contrary to the discrete-time derivation in \ref{sec:edmd}. This is done simply because of the continuous-time definition of \eqref{eq:basis}, the system can be discretized after the whole procedure.

The idea is to select a set of initial conditions $\Gamma$ for the system \eqref{eq:koop2} and simulate the system for a time $T$ with a sampling period $T_s$. Then for each sampled data point $x_{k T_s}$, evaluate a set of eigenfunctions
\begin{equation}
\psi(x_{k T_s})_{\lambda,g} = e^{k T_s} g(x_0)
\end{equation}
for some $\lambda \in \Lambda$, $g \in \mathcal{G}$, where $x_0 \in \Gamma$ is the initial condition of the state $x_{k T_s}$.

The values of $\psi(x)_{\lambda,g}$ can be interpolated in order to approximate their values for states $x$ which are not in the learning dataset.
The approximation will be denoted as $\hat{\psi}(x)_{\lambda,g}$

The set $\Gamma$ should be a non-recurrent set, meaning that the sampled trajectories of the system \eqref{eq:koop2} should not return to the states from set $\Gamma$. This condition is not difficult to fulfill with a dissipative system such a car. If, for example, the set consisted of states with the same kinetic energy, none of the trajectories would return to the set $\Gamma$, because the system naturally loses its kinetic energy as time progresses.
For more information and proofs of the above stated facts, see \cite{Korda_eigen}.

The matrix $C$ can be obtained in a similar fashion as in the previous case. To approximate the output $y = g(x)$, solve
\begin{equation}
\min_C \sum_{i=1}^{M} || g(x_i) - C\boldsymbol{\hat{\psi}}(x_i) ||_2^2
\label{eq:app_error}
\end{equation}
where $M$ is the total number of sampled states.

The controlled case is not considered here, for more information on the derivation of the $B$ matrix, see \cite{Korda_eigen}.

\subsubsection*{Results}
The simulation experiments were performed with the model described in \ref{sec:3state}. 
The non-recurrent set $\Gamma$ was chosen as set of states with the same kinetic energy, the reasoning for this selection was mentioned above.

\begin{figure}[tb]
\centerline{\includegraphics[width=0.5\textwidth]{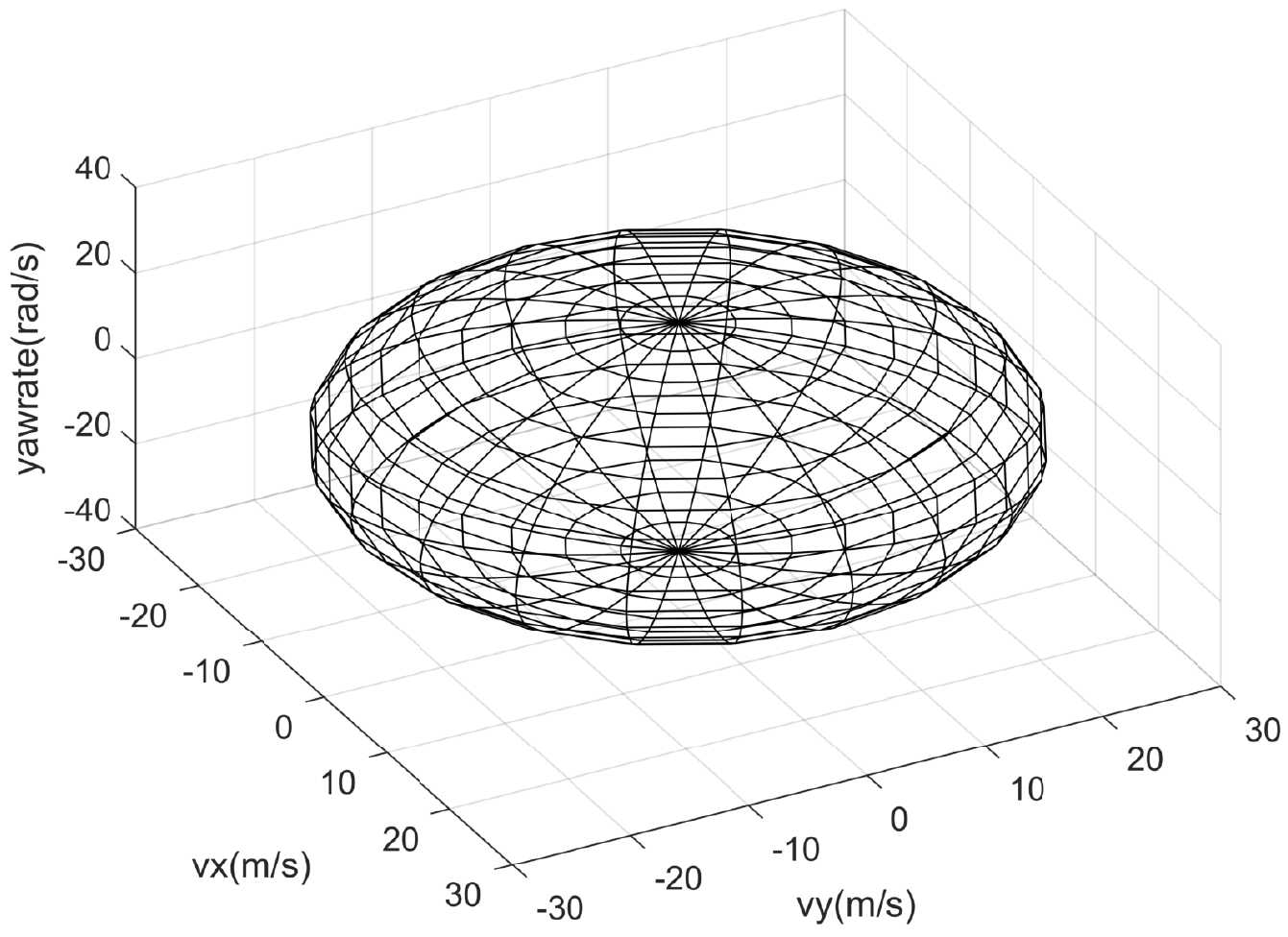}}
\caption{Set of states with constant kinetic energy equivalent to $1300kg$ car riding at $100km/h$. Each vertex of the mesh is one of the 441 initial states from which the learning dataset was constructed.}
\label{fig:kine}
\end{figure}

The level of the kinetic energy was set to $500kJ$, which is an energy equivalent of a car weighting $1300kg$, riding at $100km/h$. The set is depicted in Fig. \ref{fig:kine}.

Next, the nonlinear system was simulated from 441 initial states from the set $\Gamma$ (the set $\Gamma$ is visualized in Fig.\ref{fig:kine}) for sufficiently long time $T$ in order to cover the state space with enough data, in our case $T =1s$. The trajectories were sampled at sampling rate $T_s = 0.01s$.

Then the eigenfunctions $\psi_{\lambda,g}$ were calculated. The eigenvalues were chosen in a similar fashion as in \cite{Korda_eigen}. By applying Dynamic Mode Decomposition (DMD) algorithm \cite{dmd} to the dataset, three eigenvalues $\Lambda_{DMD}$ were obtained. Additional eigenvalues were constructed as linear combinations of elements from $\Lambda_{DMD}$, more in \cite{Korda_eigen}.

For the functions $g$, thin plate spline radial basis functions were used:
\begin{equation}
g(x) = || x - x_c ||^2 \log(|| x - x_c ||)
\end{equation}
The centers $x_c$ were chosen randomly with normal distribution $N(\mu,\sigma^2)$, where $\mu$ and $\sigma$ are the mean and standard deviation of the whole dataset.
Note that instead of interpolating the functions $\psi_{\lambda,g}$, nearest neighbour was chosen. This was done for simplicity and it yielded good results.

The approximated system was evaluated using 2000 randomly selected initial conditions, all within the ellipsoid in Fig. \ref{fig:kine}. Both the Koopman system and the real one were then simulated for $0.5s$ (which is unnecessarily long for MPC control, but it demonstrates the capabilities of the approach). The error of each trajectory is again calculated using RMSE, defined in \eqref{eq:rmse}.

\begin{figure}[t]
\centerline{\includegraphics[width=0.5\textwidth]{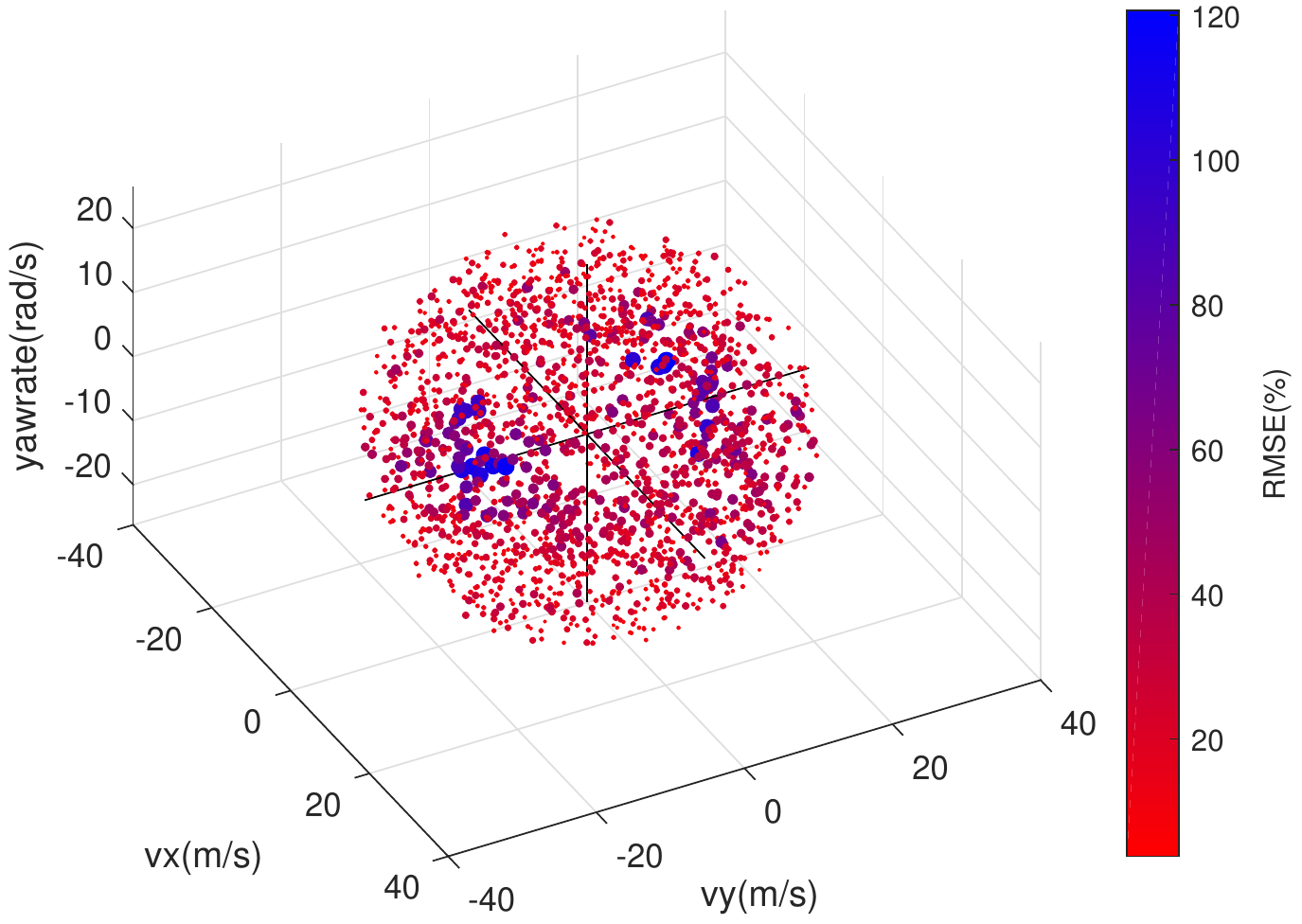}}
\caption{Each ball in the figure is an initial condition. Its size and color indicates prediction error of its associated trajectory. The mean of RMSE is $23\%$, standard deviation $15\%$.}
\label{fig:koop2_top}
\end{figure}

\begin{figure}[t]
\centerline{\includegraphics[width=0.5\textwidth]{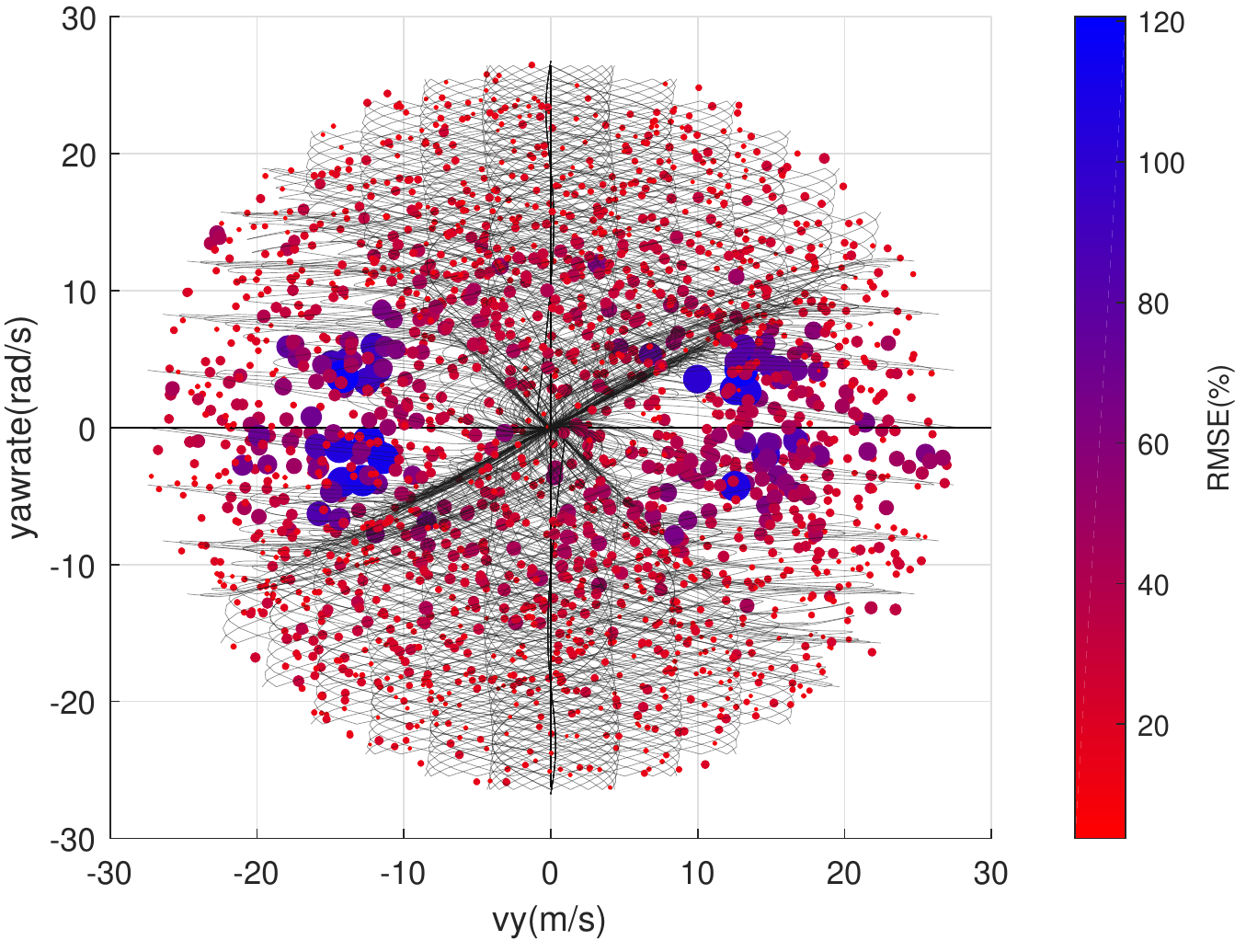}}
\caption{Trajectories used for construction of eigenfunctions are depicted as gray lines. The prediction error is large in areas with very little data. The trajectories converge to planes, forming X-shape around origin. The shape is caused by the vehicle not having enough energy to spin and going into a drift while decreasing $v_y$ and $\dot{\psi}$.}
\label{fig:koop2_side}
\end{figure}

\begin{figure}[htbp]
\centerline{\includegraphics[width=0.5\textwidth]{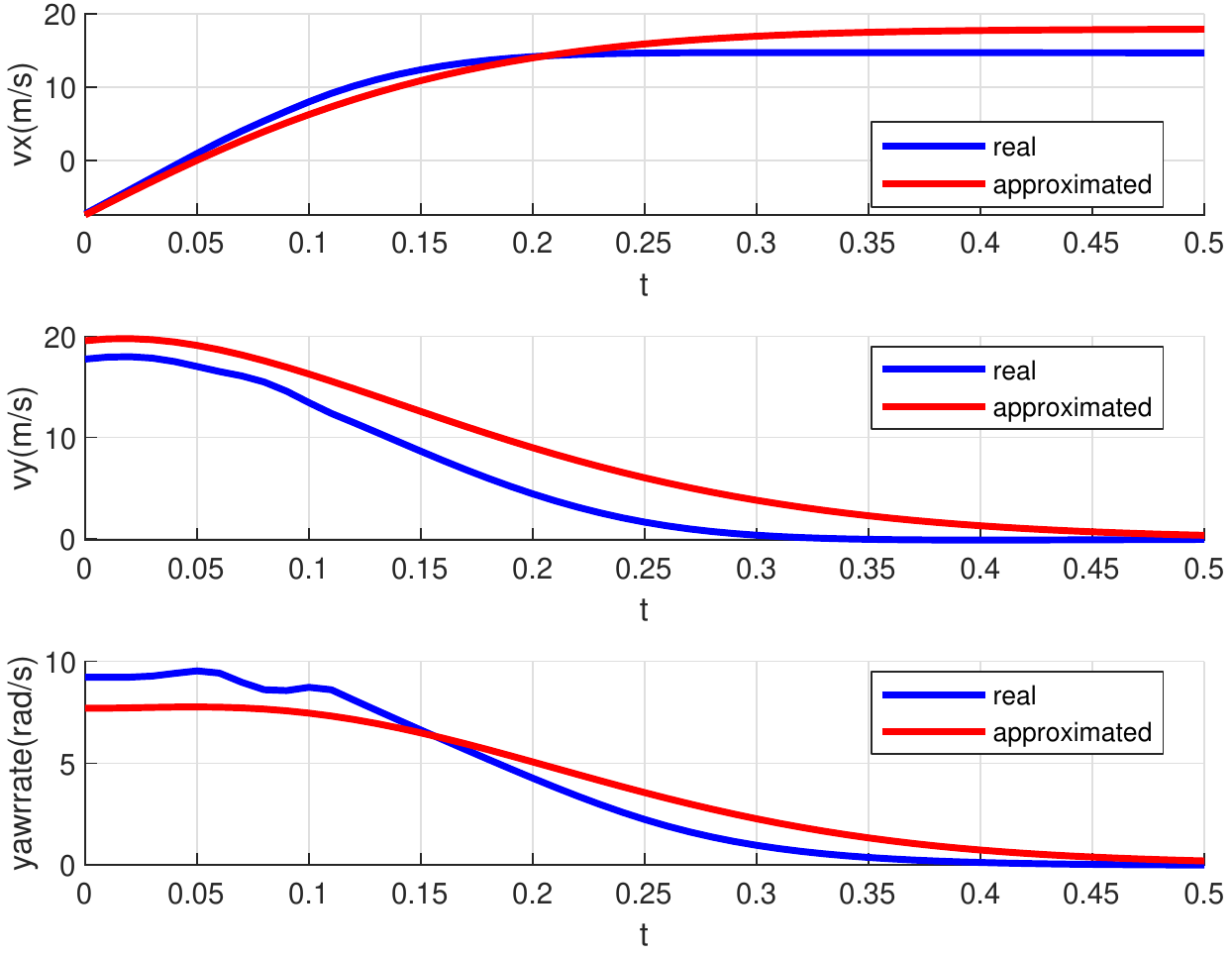}}
\caption{This figure shows a trajectory with $\text{RMSE} = 23\%$ which is equal to the mean RMSE of the whole testing dataset.}
\label{fig:koop2_mean_ex}
\end{figure}

Results can be seen in Fig.\ref{fig:koop2_top}. There are two clusters of initial conditions with large prediction error. These are the areas with very little data as can be seen in Fig.\ref{fig:koop2_side}. Trajectory with RMSE equal to the mean RMSE can be seen in Fig.\ref{fig:koop2_mean_ex}.

Notice that that in Fig.\ref{fig:koop2_mean_ex}, the initial conditions are different. This is caused by the approximation error of $\boldsymbol{\hat{\psi}}$ which was minimized in \eqref{eq:app_error}.

\section{Conclusion}
This paper presented two algorithms based on the Koopman operator for identifying vehicle dynamics.

The EDMD approach in \ref{sec:edmd} showed good result, but only with the model \ref{sec:simplified} which does not include the tire nonlinearities which are the greatest challenge in vehicle dynamics.

The eigenfunction approach in \ref{sec:eigenfun} yielded very promising results on a model with the tire nonlinearities \ref{sec:3state}.
This approach will be the focus of our future work including the controlled case, which wasn't discussed in this work, and its implementation with linear MPC.

\section*{Acknowledgment}
This research was supported by the Czech Science Foundation (GACR) under contract No. 19-16772S.
We would like to thank Milan Korda for valuable discussions on the topic of eigenfunctions and Koopman theory.

\bibliographystyle{IEEEtran}
\bibliography{../../ref} 


\end{document}